\documentclass[11pt,a4paper]{amsart}
\usepackage[margin=2.5cm]{geometry}
\usepackage{hyperref,enumerate,xcolor} 
\usepackage{amscd, amssymb,braket}
\usepackage{tikz-cd}
\usepackage{graphicx}
\usepackage[all]{xy}

\usepackage[nameinlink]{cleveref}
\usepackage{verbatim}

\newtheorem{thm}[equation]{Theorem}
\newtheorem{prop}[equation]{Proposition}

\newtheorem{definition}[equation]{Definition}
\newtheorem{lemma}[equation]{Lemma}

\numberwithin{equation}{section}

\newcommand{\A}{{\mathbb A}}
\newcommand{\C}{{\mathbb C}}

\newcommand{\Spin}{\mathrm{Spin}}

\newcommand{\PGSp}{\mathrm{PGSp_6}}

\newcommand{\PGSO}{\mathrm{PGSO}}

\newcommand{\SO}{\mathrm{SO}}

\newcommand{\PD}{\mathrm{PD}}
\newcommand{\Sp}{\mathrm{Sp_6}}
\newcommand{\PGL}{\mathrm{PGL}}
\newcommand{\GL}{\mathrm{GL}}
\newcommand{\SL}{\mathrm{SL}}

\newcommand{\Aut}{\mathrm{Aut}}

\newcommand{\Ind}{\mathrm{Ind}}
\newcommand{\Irr}{\mathrm{Irr}}

\title{A Theory of $\gamma$-factors  for $G_2 \times \GL_r$}
 
\author{Wee Teck Gan and Gordan Savin}
\address{W.T.G.:   Department of Mathematics, National University of Singapore, 10 Lower Kent Ridge Road
Singapore 119076} \email{matgwt@nus.edu.sg}
\address{G. S.: Department of Mathematics, University of Utah, Salt Lake City, UT 84000, USA}\email{savin@math.utah.edu}

\begin{document}
\maketitle
\section{\bf Introduction}
As the title of this paper indicates, our purpose is to develop a theory of local $\gamma$-factors
\[ \gamma(s, \pi \times \tau, \psi)  \]
 for irreducible representations  $\pi \otimes \tau$ of the group $G_2 \times \GL_r$ over a local field $F$, associated to the tensor product of the 7-dimensional representation of the Langlands dual group $G_2^{\vee} = G_2(\C)$ and the $r$-dimensional standard representation of  $\GL_r^{\vee} = \GL_r(\C)$. The reader may be expecting that this involves the discovery of  an ingenious  global  zeta integral (in the sense of Tate's thesis) which represents the relevant L-function and the proof  of a local functional equation for the corresponding local zeta integral which can  then be used to define the desired $\gamma$-factor. Such a reader will be disappointed, for this is not the route we take. Rather, we define these $\gamma$-factors by establishing a functorial lifting from $G_2$ to the group $\GL_7$, associated to the embedding 
\[  G_2(\C) \longrightarrow \GL_7(\C) \]
of dual groups, and then defining the $\gamma$-factor for $G_2 \times \GL_r$ as the Rankin-Selberg $\gamma$-factor of $\GL_7 \times \GL_r$.  The functorial  lifting needed for this definition is  provided by   the exceptional theta correspondences that we established in \cite{GS}, together with the local lifting from $\Sp$ to $\GL_7$ (which is a special case of results of Arthur \cite{Ar}). 
\vskip 5pt

Now the knowledgeable reader will no doubt question the utility of such an approach, so let us elaborate further. Any such potential theory of $\gamma$-factors is expected to satisfy a list of properties \cite{Sh, LR}. The most important of these properties can be loosely summarized as:
\vskip 5pt

\begin{itemize}
\item {\em Unramified property}: the $\gamma$-factors for unramified representations must be the expected ones, in the sense that they should agree with what is provided by the unramified local Langlands correspondence.
\vskip 5pt
\item {\em Multiplicativity}:  the $\gamma$-factors should satisfy an inductive property with respect to parabolic induction. In particular, the $\gamma$-factors depend only on the cuspidal support of the representations and the $\gamma$-factor of non-supercuspidal representations can be deduced from simpler $\gamma$-factors of supercuspidal representations on proper Levi subgroups.

\vskip 5pt

\item {\em Global property}: over a number field $k$, the $\gamma$-factors can be used to provide the global functional equation for the corresponding global L-function of cuspidal representations. 
\end{itemize}
Equally importantly, the above properties uniquely characterize the relevant $\gamma$-factors. 
\vskip 5pt

Given this, one would naturally ask if the $\gamma$-factors for $G_2 \times \GL_r$ whose definition we outlined above satisfy such a list of properties and is uniquely characterized by it. 
This is in fact the main result of this paper. More precisely,  we show in Theorem \ref{T:gamma} that the $\gamma$-factors for $G_2 \times \GL_r$ that we defined satisfy a list of properties as above, with the caveat that the global property above is only established for a subset of cuspidal representations, instead of all of them. This may seem a less-than-ideal result, but Theorem \ref{T:gamma} further asserts that this weaker list of properties already characrerizes the local $\gamma$-factors. Hence, if there is indeed  a theory of $\gamma$-factors for $G_2 \times \GL_r$, the one we developed here is the only possible candidate.  Whether it satisfies the full list of expected properties (for all cuspidal representations) remains to be seen. 
\vskip 5pt

While Theorem \ref{T:gamma} provides a satisfactory justification for our definition of $\gamma$-factors, one may ask whether such a construction lends itself readily to concrete applications. 
If one is thinking about arithmetic applications such as relating L-values to period integrals, then our construction indeed has very little to offer in that direction. However, the main motivation for constructing this theory of $\gamma$-factors is its use in characterizing the local Langlands correspondence (LLC) for $G_2$ that we established in \cite{GaSa2}.
See \cite[Main Theorem (ix) and (x)]{GaSa2}.
\vskip 5pt

Using the same strategy, we also define the twisted adjoint $\gamma$-factor 
\[   \gamma(s, \pi, {\rm Ad} \times \chi, \psi) \]
for $ \pi \otimes \chi \in \Irr(G_2 \times \GL_1)$ and show that they are characterized by an analogous list of properties. We also deduce that the LLC for $G_2$ respects these $\gamma$-factors as well.  
\vskip 5pt

\vskip 10pt

 \noindent{\bf Acknowledgments:} 
W.T. Gan is partially supported by a Singapore government MOE Tier 1 grant R-146-000-320-114 and a Tan Chin Tuan Centennial Professorship.  
G. Savin is partially supported by a Gift No. 946504 from Simons Foundation. 
\vskip 10pt

\section{\bf $\gamma$-Factors for $G_2 \times \GL_r$}
In this section, we give the full details of our construction of a theory of $\gamma$-factors for $G_2 \times \GL_r$ outlined in the introduction.
\vskip 5pt

\subsection{\bf Recollection of \cite{GS}}
We first recall some results from \cite{GS} over a nonarchimedean local field $F$.
The main results of \cite{GS} can be summarized as the construction of a map
\[  {\rm Lif} : \Irr(G_2) \longrightarrow  \Irr(\GL_7) \]
which is given in the following commutative diagram:

\vskip 5pt
\[
\begin{CD}
   {}   @. \Irr(G_2)   @. {}  \\
   @.  @|  @. \\
\Irr(PD^{\times}) @<\theta<<  \Irr^{\spadesuit}(G_2) \sqcup \Irr^{\heartsuit}(G_2) @>\theta>> \Irr(\PGSp) \\
@VV{\rm JL}V     @VV{\rm Lif}V  @VV{\rm rest}V  \\
\Irr_{ds}(\PGL_3) @>{\boxplus}>>  \Irr(\GL_7)  @<{\rm Art}<<   \Irr(\Sp)/_{\PGSp}  
\end{CD} \]
Here, 
\vskip 5pt
\begin{itemize}
\item the two $\theta$'s refer to the local theta correspondence;
\item $D$ is a cubic division $F$-algebra, so that $PD^{\times}$ is an inner form of $\PGL_3$;
\item $ \Irr^{\spadesuit}(G_2)$ and $\Irr^{\heartsuit}(G_2)$ are the subsets of $\Irr(G_2)$ which participate in the theta correspondence with $\PD^{\times}$ and $\PGSp$ respectively.
\item JL refers to the Jacquet-Langlands transfer from $PD^{\times}$ to $\PGL_3$;
\item rest refers to the restriction of representations of $\PGSp$ to $\Sp$;
\item Art refers to the functorial transfer of irreducible representations from $\Sp$ to $\GL_7$ established by Arthur \cite{Ar};
\item $\boxplus: \Irr_{ds}(\PGL_3) \longrightarrow \Irr(\GL_7)$ is the map 
\[  \tau \mapsto \tau \times 1 \times \tau^{\vee}. \]
\end{itemize}
\vskip 5pt

\subsection{\bf Definition of $\gamma$-factors}
With the above preparation, we can now define the local $\gamma$-factors.
\vskip 5pt

\begin{definition}  \label{D:gamma}
Fix a nontrivial additive character $\psi: F \rightarrow \C^{\times}$. 
For $\pi \in {\rm Irr}(G_2)$ and $\rho \in {\rm Irr}(\GL_r)$, set
\[  \gamma(s, \pi \times \rho,\psi) :=  \gamma(s, {\rm Lif}(\pi) \times \rho, \psi). \]
where the $\gamma$-factor on the RHS is the Rankin-Selberg $\gamma$-factor defined by Jacquet-PS-Shalika and Shahidi \cite{GeSh}.
\end{definition}
\vskip 5pt

\subsection{\bf Compatibility with LLC}  In \cite[Main Theorem (ix)  and \S 3.4]{GaSa2}, we have shown that the LLC for $G_2$ established there respects the above $\gamma$-factor:
\vskip 5pt

\begin{thm}  \label{T:comp-LLC}
Let $F$ be a non-archimedean local field. 
For $\pi \in \Irr(G_2)$ and  $\rho \in \Irr(\GL_r)$, one has
\[  \gamma(s, \pi \times \rho, \psi) = \gamma(s,  ({\rm std} \circ \phi_{\pi}) \otimes \phi_{\rho}, \psi) \]
 where $\phi_{\pi}: WD_F \longrightarrow G_2(\C)$ and $\phi_{\rho}: WD_F \longrightarrow \GL_r(\C)$  are the L-parameters of $\pi$ and $\rho$ respectively and ${\rm std}: G_2(\C) \longrightarrow \GL_7(\C)$ is the standard degree 7 irreducible  representation of $G_2(\C)$. 
\end{thm}

  \vskip 5pt

The reader will no doubt complain that our construction of the LLC map in \cite{GaSa2} and our definition of the $\gamma$-factor here makes the above theorem almost a formality, which would not be inaccurate description. But the main question is then whether this $\gamma$-factor can be characterized by a list of properties independent of the LLC.
 Our main result, Theorem \ref{T:gamma} below, gives an affirmative answer. 
Before coming to that, we  define $ \gamma(s, \pi \times \rho,\psi)$ in the archimedean setting for the sake of completeness.
\vskip 5pt

\subsection{\bf Archimedean case}
 In the archimedean setting, the (known) LLC allows one to define 
\[  {\rm Lif} : \Irr(G_2) \longrightarrow \Irr(\GL_7). \]
More precisely, if $\phi_{\pi}$  is the L-parameter of $\pi$, then ${\rm Lif}(\pi)$ is characterized by
\[  \phi_{{\rm Lif}(\pi)}  = {\rm std} \circ \phi_{\pi}. \]
 If $\phi_{\rho}$ denotes  the L-parameters of $\rho \in \Irr(\GL_r)$, then we set:
 \[  \gamma(s, \pi \times \rho,\psi) :=\gamma(s, {\rm Lif}( \pi) \times \rho, \psi) = \gamma(s, \phi_{\pi} \otimes \phi_{\rho},\psi). \]
where  the second $\gamma$-factor is that defined by Jacquet-PS-Shalika and Shahidi, and for the third $\gamma$-factor, we have regarded  $\phi_{\pi}$ and $\phi_{\rho}$ as $7$-dimensional and $r$-dimensional representations respectively.
\vskip 5pt
 
 The main difference with the non-archimedean case is  that we have used the LLC instead of the theta correspondence to define the lifting map ${\rm Lif}$.
 
 \vskip 5pt
 
 \subsection{\bf Cuspidal support} Returning to the setting of a non-archimedean local field $F$, we shall now investigate the behaviour of cuspidal support under the map ${\rm Lif}$.
  \vskip 5pt
  
Let $\mathcal{SC}(G)$ denote the set of cuspidal supports of non-supercuspidal representations for the group $G$:
\[  \mathcal{SC}(G) = \{ [ M, \tau]: \text{$M$ is a proper Levi subgroup of $G$ and $\tau \in \Irr_{sc}(M)$} \} \]
where the pairs $[M,\tau]$ are taken up to $G(F)$-conjugacy. 
\vskip 5pt

For the group $G_2$, the cuspidal supports fall into 3 families, as there are 3 conjugacy classes of parabolic subgroups:
\vskip 5pt
\begin{itemize}
\item $[T, \chi]$ with $T$ a maximal torus;
\item $[M, \tau]$, with $M \cong \GL_2$ the Levi subgroup of a Heisenberg maximal parabolic subgroup $P$;
\item $[L,\tau]$, with $L \cong \GL_2$ the Levi subgroup of a non-Heisenberg maximal parabolic subgroup $Q$.
\end{itemize}
Then we define
\[  f:  \mathcal{SC}(G_2) \longrightarrow \mathcal{SC}(\GL_7)  \]
as follows:
\begin{itemize}
\item $f(T, \chi)  = (T' ,\chi')$, where $T'$ is the maximal diagonal torus of $\GL_7$ and $\chi'$ is given as follows. Let us identify the character $\chi$ with a triple 
$(\chi_1, \chi_2, \chi_3)$ of characters of $F^{\times}$ with $\chi_1\chi_2\chi_3 =1$, as in \cite[Pg. 16]{GS}. Then
\[  \chi' =  \chi_1  \times \chi_2\times \chi_3 \times 1 \times \chi_3^{-1} \times \chi_2^{-1} \times \chi_1^{-1}. \]

\item 
\[ f( M, \tau) = (\GL_2 \times \GL_2 \times \GL_1^3, \, \tau \otimes \tau^{\vee} \otimes \omega_{\tau} \otimes \omega_{\tau}^{-1} \otimes 1) \]

\item 
\[  f(L,\tau)  = (\GL_2 \times \GL_2 \times \mathcal{L}, \,  \tau \otimes \tau^{\vee} \otimes \rho_{\tau}), \]
where $(\mathcal{L}, \rho_{\tau})$ is the cuspidal support of the representation ${\rm Ad}(\tau)$ of $\GL_3$. 
If $\tau$ is a non-dihedral supercuspidal representation, then the adjoint lift ${\rm Ad}(\tau)$ is a supercuspidal representation of $\GL_3$. On the other hand, if $\tau$ is dihedral, then ${\rm Ad}(\tau)$ is non-supercuspidal  and its 
cuspidal support can be a maximal proper Levi subgroup subgroup of $\GL_3$ or even the maximal torus, depending on whether $\tau$ is dihedral with respect to 1 or 3 quadratic fields. 
\end{itemize}
Having introduced the map $f$, the following result follows from the explicit theta correspondence established in \cite{GS}:
\vskip 5pt

\begin{prop} \label{P:SC}
For any non-supercuspidal $\pi \in \Irr(G_2)$, one has
\[  \mathcal{SC}({\rm Lif}(\pi)) = f(\mathcal{SC}(\pi)). \]
In particular, if $\pi$ and $\pi'$ have the same supercuspidal support, then so do ${\rm Lif}(\pi)$ and ${\rm Lif}(\pi')$. 
\end{prop}
\vskip 5pt

\vskip 5pt

\subsection{\bf Properties and characterization of $\gamma$-factor}
The main result of this paper is the following theorem:
\vskip 5pt

\begin{thm}  \label{T:gamma}
The local $\gamma$-factors defined above for non-archimedean local fields satisfy and are uniquely characterized by the following properties:
\vskip 5pt

\noindent (i)  $\gamma(s, \pi \times \rho,\psi)$ depends only on the cuspidal support of $\pi$ and $\rho$. More precisely, 
\vskip 5pt
\begin{itemize}
\item[(a)] If $\pi$ is a subquotient of $\Ind_B^{G_2} \chi$, with $\chi = (\chi_1, \chi_2, \chi_3)$ and  $\chi_1 \chi_2 \chi_3 =1$, then
\[  \gamma(s, \pi \times \rho,\psi) = 
\gamma(s, \rho, \psi) \cdot  \prod_{i=1}^3 \left( \gamma(s,  \chi_i \times \rho) \cdot \gamma(s, \chi_i^{-1} \times \rho) \right). \]
\vskip 5pt

\item[(b)] If $\pi$ is a subquotient of $\Ind_P^{G_2} \tau$, then
\[   \gamma(s, \pi \times \rho,\psi) = \]
\[ \gamma(s, \tau \times \rho,\psi) \cdot \gamma(s, \tau^{\vee} \times \rho, \psi) \cdot \gamma(s, \omega_{\tau} \times \rho, \psi) \cdot 
\gamma(s, \omega_{\tau}^{-1} \times \rho, \psi) \cdot \gamma(s, \rho, \psi). \]
\vskip 5pt

\item[(c)] If $\pi$ is a subquotient of $\Ind_Q^{G_2} \tau$, then
\[  \gamma(s, \pi \times \rho,\psi) = \gamma(s, \tau \times \rho,\psi) \cdot \gamma(s, \tau^{\vee} \times \rho, \psi) \cdot \gamma(s, \rho_{\tau} \times \rho, \psi). \]
\vskip 5pt

\item[(d)] If $\rho$ is a subquotient of $\rho_1 \times \rho_2$, then
\[  \gamma(s, \pi \times \rho,\psi) = \gamma(s, \pi \times \rho_1,\psi) \cdot \gamma(s, \pi \times \rho_2,\psi). \]
\end{itemize}
\vskip 5pt
\noindent In particular, the $\gamma$-factors for unramified representations are compatible with the unramified LLC by (a). Moreover, the $\gamma$-factors satisfy the property of multiplicativity.  
  \vskip 5pt
  
\noindent (ii) Suppose we have the following data:
\vskip 5pt

\begin{itemize}
\item $k$ is a totally real number field with adele ring $\A$;
\item $\psi: k \backslash \A \longrightarrow \C^{\times}$ is a nontrivial additive character;
\item $\mathbb{O}$ is an (not necessarily split) octonion $k$-algebra with automorphism group $G= \Aut(\mathbb{O})$, which is of type $G_2$.
 \item $\Pi$ is a cuspidal representation of $G$ such $\Pi_v$ is spherical for all real places $v$ where $G(k_v)$ is split.  
 \item $\Sigma$ is a cuspidal representation of $\GL_r$;
\item $S$ is a finite set of places $v$ of $k$ containing all archimedean ones and such that  $\Pi_v$, $\Sigma_v$ and $\psi_v$ are unramified for $v \notin S$.
\end{itemize}
\vskip 5pt

Suppose further that one of the two conditions holds: 
\vskip 5pt
\begin{itemize}
\item[(a)]  $\Pi$ has a nonzero cuspidal global theta lift to the split group $\PGSp$. 

\item[(b)] $\Pi$ has a nonzero theta lift to 
$\mathbb{PD}^{\times}$ for some central simple $k$-algebra $\mathbb{D}$  of degree $3$ (possibly split).
\end{itemize} 

\vskip 5pt

Then $L^S(s, \Pi \times \Sigma)$ has meromorphic continuation to $\C$ and satisfies the functional equation
\[ L^S(1-s, \Pi \times \Sigma^{\vee}) =  \left( \prod_{v \in S}    \gamma(s, \Pi_v \times \Sigma_v,\psi_v) \right) \cdot 
L^S(s, \Pi \times \Sigma). \]
 
\end{thm}
\vskip 5pt
Before proving the theorem, we make two remarks:
\vskip 5pt

\noindent{\bf Remarks:} 
\vskip 5pt
\noindent (i)  It is curious to note that Theorem \ref{T:gamma} does not include the archimedean case.  The main reason is that we are not entirely sure if 
  the key property  of multiplicativity is satisfied. Namely if $\pi_1$ and $\pi_2$ are constituents of the same principal series representation,  does one have
 \[   \gamma(s, \phi_{\pi_1} \otimes \phi_{\rho},\psi) =  \gamma(s, \phi_{\pi_2} \otimes \phi_{\rho},\psi)? \]
 The two papers \cite{A1, A2}  of Martin Andler come very close to giving an affirmative answer to this question, but seem to be  not as definitive as one might hope. This is an issue which certainly deserves to be definitively resolved, but we will not discuss it further here.
\vskip 10pt

\noindent (ii) The properties stated in Theorem \ref{T:gamma} are the usual ones used for characterizing $\gamma$-factors,  except that the global property (ii) is not the most general. One would like to have the global property (ii) for all cuspidal representations $\Pi$ of $G_2$. 
One reason for this more restricted result is due to the fact that the relevant archimedean exceptional theta correspondences are not fully understood. Another is that the global lifting from $G_2$ to $\GL_7$ cannot be fully achieved by the exceptional theta correspondences considered here.
On the other hand, the uniqueness part of Theorem \ref{T:gamma}, proved under this more constrained  global property,  is thus stronger than the usual version and implies that the theory of $\gamma$-factors developed here is the only possible candidate for such a theory.
\vskip 10pt

\subsection{\bf Proof of Expected Properties}
We now show that the $\gamma$-factors defined in Definition \ref{D:gamma} satisfy properties (i) and (ii) of the theorem.
The properties in (i)  follow from the definition of the $\gamma$-factors, by using Proposition \ref{P:SC} and the multiplicativity of the Rankin-Selberg $\gamma$-factors. 
\vskip 5pt

Let us now place ourselves in the setting of (ii), so that $\Pi$ is a cuspidal representation of $G= \Aut(\mathbb{O})$ over $k$. We consider the situations (a) and (b) in turn.

\vskip 5pt
\begin{itemize}
\item[(a)]  Suppose first that we are in the situation (a)  in (ii) of the theorem, so that 
$\Pi$ has nonzero cuspidal global theta lift $\Theta(\Pi)$ on $\PGSp$.  
\vskip 5pt

The local theta correspondence at the real places $v$ such that $G(k_v)$ is compact is completely understood by the work of Gross-Savin \cite{GS97}. 
The local theta lift of an irreducible finite-dimensional representation of compact $G_2$ is the holomorphic discrete series of $\PGSp(\mathbb{R})$ with the associated  infinitesimal character.   If $G(k_v)$ is split and $\Pi_v$ is spherical then $\Theta(\Pi_v)$ is a finite length module with unique (spherical) quotient by Loke-Savin \cite{LS}.  In these two cases 
the local theta correspondence is functorial with respect to the inclusion $G_2(\mathbb{C}) \hookrightarrow \Spin_7(\mathbb{C})$. 
\vskip 5pt

In view of the above, $\Pi'=\Theta(\Pi)$ is 
irreducible and $\Pi'_v \cong \theta(\Pi_v)$ at all places. 
Let  $\mathcal{A}(\Pi' )$ denote its Arthur transfer to $\GL_7$.  
 \noindent Thus, we see that at all places $v$ of $k$, 
 \[  \mathcal{A}(\Pi')_v \cong  {\rm Lif}(\Pi_v). \]
 In particular, 
 \[  \gamma(s, \Pi_v \times \Sigma_v, \psi_v)= \gamma(s, \mathcal{A}(\Pi'_v) \times \Sigma_v, \psi_v) \]
 and
 \[  L^S(s, \Pi \times \Sigma) =  L^S(s, \mathcal{A}(\Pi'_v)\times \Sigma). \]
 so the desired global functional equation follows from that for the Rankin-Selberg L-functions for $\mathcal{A}(\Pi') \otimes \Sigma$ on $\GL_7 \times \GL_r$. 
   \vskip 5pt

\item[(b)]  Suppose now that we are in situation (b) in (ii) of the theorem.  
In this case $\Pi$ has a global theta lift $\Theta(\Pi)$ on $\Aut(D)$.  
\vskip 5pt

At real places, the Howe duality theorem holds for spherical representations by results of Loke-Savin \cite{LS}. In particular, $\Theta(\Pi)$ is irreducible and all its 
local components are isomorphic to $\theta(\Pi_v)$.  At any real place $v$, 
 the L-parameters of  $\theta(\Pi_v)$ and $\Pi_v$ are related by the inclusion  
 $\SL_3(\mathbb{C})\rtimes \mathbb Z/2\mathbb Z \hookrightarrow G_2(\mathbb{C})$ of dual groups. 
 \vskip 5pt
 
 The lift of $\Pi$ to $\mathbb P\mathbb D^{\times}$  is obtained by restricting $\Theta(\Pi)$ to $\mathbb P\mathbb D^{\times}$. Let $\Pi'$ be any irreducible summand. 
Then $\Pi'$ has  a  Jacquet-Langlands transfer $\mathrm{JL}(\Pi')$ to $\PGL_3$.  At all places $v$, we thus have
 \[
   \gamma(s, \Pi_v \times \Sigma_v, \psi_v) =\gamma(s, \psi_v) \cdot  \gamma(s, \mathrm{JL}(\Pi')_v \times \Sigma_v, \psi_v)  \cdot \gamma(s, \mathrm{JL}(\Pi')^{\vee}_v  \times \Sigma_v, \psi_v).  \]
Now the desired global functional equation for $\Pi \times \Sigma$ follows from that for the Rankin-Selberg L-function of  $JL(\Pi') \times \Sigma$ on  $\GL_3 \times \GL_r$.

 \end{itemize}
The completes the proof that the $\gamma$-factors  defined in Definition \ref{D:gamma} satisfy properties (i) and (ii) of the theorem.
\vskip 10pt

\subsection{\bf Proof of Uniqueness}
Suppose that we have two systems of invariants $\gamma_1$ and $\gamma_2$ as in the theorem. 
By properties (a), (b) and  (c) in (i), we   deduce the equality  
\[   \gamma_1(s, \pi \times \rho,\psi) = \gamma_2(s, \pi \times \rho,\psi)  \]
when $\pi$ is non-supercuspidal. When $\pi$ is supercuspidal, property (d) of (i) reduces us to showing the above equality when 
$\rho$ is supercuspidal as well. For the supercuspidal case, we shall exploit the global property (ii). Let $k$ be a totally real field such that $k_{v_0}\cong F$ for some  place $v_0$.
Let $\Sigma$ be a cuspidal representation of $\GL_n(k)$ such that $\Sigma_{v_0}\cong \rho$. 
There are now  two cases:

\vskip 5pt
\noindent \underline{Case (i)}: $\pi \in   \Irr_{sc}^{\heartsuit}(G_2(F))$. Thus the theta lift $\theta(\pi)$ to $\PGSp$ is non-zero. It is a tempered representation, since this 
theta correspondence preserves tempered representations. Let $P=MN$ be a Siegel parabolic in $\PGSp$. 
Open $M$-orbits on $\hat N$ are parameterized by quaternion algebras.
Since $\theta(\pi)$ is tempered, there exists a quaternion algebra $B$ such that 
$\theta(\pi)_{N,\psi_B}\neq 0$, where $\psi_B$ is a character of $N$ in the $M$-orbit parameterized by $B$. 
(Indeed, non-trivial representations of $\Sp$  such that twisted $N$-coinvaraints vanish for characters in open orbits are theta lifts from $\mathrm{O}(2)$ \cite{Li}, and these are not tempered.) 
Let $O$ be an octonion algebra over $F$, so that 
$\Aut(O)$ is the split group of type $G_2$. Then $B\hookrightarrow O$, embedding unique up to conjugation. Then $\pi$ has a non-trivial  
$\Aut(B\hookrightarrow O) \cong SL_1(B)$-invaraint linear form.  Pick a totally definite octonion algebra $\mathbb O$ over $k$, that is, $\mathbb O_v$ is definite for all real places $v$. 
Pick a totally definite quaternion algebra $\mathbb B$ such that $\mathbb B_{v_0}\cong B$. Then $\mathbb B\hookrightarrow \mathbb O$. 
By a result of Prasad-Schulze-Pillot \cite{PSP}, one can find an 
automorphic representation  $\Pi$ of $G=\Aut(\mathbb O)$ such that  
\vskip 5pt 
\begin{itemize}
\item  $\Pi_{v_0} \cong \pi$;
\item  $\Pi_v$ is a spherical representation for all finite $v \ne v_0$;
\item  $\Pi$ has a  nonvanishing global $\Aut(\mathbb B\hookrightarrow \mathbb O)$-period. 
 \end{itemize}
\vskip 5pt 
Now $\Pi$ has a non-zero lift to $\PGSp$, since it has a non-vanishing $\psi_{\mathbb B}$-Fourier coefficient, see \cite{GS97}.  Assume that $\theta(\pi)$ is supercuspidal. Then $\Theta(\Pi)$ is a cupidal automorphic representation and 
we can apply the global functional equation (assumption (ii) (a) in the theorem) 
to the two system of invariant $\gamma_1$ and $\gamma_2$, we obtain:
\[   L^S(1-s, \Pi \times \Sigma^{\vee}) =  \left( \prod_{v \in S}    \gamma_1(s, \Pi_v \times \Sigma_v,\psi_v) \right) \cdot 
L^S(s, \Pi \times \Sigma) \]
and
\[ L^S(1-s, \Pi \times \Sigma^{\vee}) =  \left( \prod_{v \in S}    \gamma_2(s, \Pi_v \times \Sigma_v,\psi_v) \right) \cdot 
L^S(s, \Pi \times \Sigma). \]
Since $\Pi_v$ is non-supercuspidal at all finite $v \ne v_0$,  the two invariants  agree at these places.  At the real places, we have taken the local gamma factors defined via the LLC, so that the two invariants are the same as well. We then  conclude from the two global functional equations above that the two $\gamma$-factors agree at the place $v_0$, i.e.
\[   \gamma_1(s, \pi \times \sigma,\psi) = \gamma_2(s, \pi \times \sigma,\psi).  \]
If $\theta(\pi)$ is not supercuspidal, then we can construct a $\Pi$ so that one other place $v_1$ the local component $\Pi_{v_1}$  is a supercuspidal representation 
that lifts to a supercuspidal representation of $\PGSp$ (and for which we already proved equality of gamma factors). 
\vskip 5pt

\vskip 5pt 
Observe that we have so far proved uniqueness of gamma factors for all $\pi \in   \Irr^{\heartsuit}(G_2(F))$, in particular for $\pi$ that are lifts from 
$\PGL_3$. This will be used in the second case: 

\vskip 5pt
\noindent \underline{Case (ii)}: $\pi \in   \Irr_{sc}^{\spadesuit}(G_2(F))$.
 Recall that $k$ is a totally real number field with a finite place $v_0$ such that $k_{v_0} \cong F$.

 \begin{lemma} 
Let $\mathbb{O}$ be the split octonion $k$-algebra. There is a cubic division algebra $\mathbb D$ ramified at $v_0$ and one another place $v_1$ and 
a cuspidal representation $\Pi$ of $G = \Aut(\mathbb{O})$ such that 
\vskip 5pt
\begin{itemize}
\item $\Pi_{v_0} \cong \pi$ and $\Pi_{v_1}$ is non-supercuspidal;
\item $\Pi_v$ is spherical for all real places $v$ of $k$;
\item the global theta lift of $\Pi$ to $\mathbb P\mathbb{D}^{\times}$ is not 0. 
\end{itemize}
\end{lemma} 

\vskip 5pt

\begin{proof}  
 We have $\pi = \theta(\tau)$  for $\tau \in \Irr(\PD^{\times})$.  Let $P=MN$ be a Heisenberg maximal parabolic in $G$. Open $M$-orbits of characters of $N$ are parameterized 
 by  \'etale cubic algebras. There exists a cubic field extension 
 $E \hookrightarrow D$ of $F$ such that $\pi_{N, \psi_{E}} \neq 0$, where $\psi_{E}$ is a character of $N$ in the $M$-orbit parameterized by $E$, \cite[Prop. 3.4 ]{GS}. 
   Then we know that $\tau$ has a nonzero 
 $PE^{\times}=\Aut(E\hookrightarrow  D)$-invariant linear form. 
 Let $e$ be a cubic field extension of $k$ such that $e_{v_0} \cong E$.  Let $v_1$ be any other place such that $e_{v_1}$ is a field. Let 
 $\mathbb D$ be a cubic division algebra ramified at $v_0$ and $v_1$.  Then $e \hookrightarrow \mathbb{D}$.   
 Let $H=\Aut(\mathbb D)$ where $\mathbb D$ is considered a Jordan algebra. Then $H(k_{v_0})\cong \PD^{\times}$ but $H(k_v)\cong \PGL_3(k_v) \rtimes \mathbb Z/2\mathbb Z$ 
for all $v\neq v_0,v_1$. 
 
 \vskip 5pt 
 
 By the globalization result of Prasad-Schulze-Pillot \cite{PSP}, one can find a cuspidal representation 
 $\Xi$ of $H$ such that 
\vskip 5pt
\begin{itemize}
\item $\Xi_{v_0} \cong \tau$ and  $\Xi_{v_1} =1$ is the trivial representation;
\item $\Xi_v$ is spherical for each real place $v$;
\item $\Xi$ has a nonvanishing global $\Aut(e\hookrightarrow \mathbb D)$-period. 
\end{itemize}
Consider now the global theta correspondence for the dual pair $H \times G$ in a split rank 2 $E_6$ constructed by means of $\mathbb D$. 
The global theta lift $\Theta(\Xi)$ of $\Xi$ to $G_2$ is thus nonzero, since it has nonvanishing $\psi_{e}$-Fourier coefficient.  
Moreover, $\Theta(\Xi)$ is contained in the space of cusp forms as well (because of the supercuspidality of $\pi$ at the place $v_0$).  
At finite places $v$ of $k$, the local theta correspondence is fully understood by \cite{GS}. On the other hand, at real places, the Howe duality theorem holds for spherical representations by results of Loke-Savin \cite{LS}. Thus  $\Pi=\Theta(\Xi)$ is an irreducible cuspidal representation of $G$ which satisfies:
\vskip 5pt

\begin{itemize}
\item $\Pi_{v_0}  \cong \pi$ and $\Pi_{v_1}$ is the local theta lift of the trivial representation of $PD^{\times}$, which is a non-supercuspidal discrete series representation; 
\item $\Pi_v$ is spherical at all real places $v$; 
\item the global theta lift of $\Pi$ to $\mathbb{PD}^{\times}$ is equal to the restriction of $\Xi$ to $\mathbb{PD}^{\times}$ and hence is not 0. 
\end{itemize}
In particular, $\Pi$ satisfies the requirements of the lemma.
\end{proof}

Since $\Pi_{v_1}$ is not supercuspidal, and $\Pi_{v} \in  \Irr^{\heartsuit}(G_2(k_v))$  for each finite place $v\neq v_0,v_1$ we have 
\[ 
\gamma_1(s, \Pi_v \times \Sigma_v,\psi_v) =\gamma_2(s, \Pi_v \times \Sigma_v,\psi_v) 
\] 
for all finite places $v\neq v_0$. Now the global functional equations for $\gamma_1$ and $\gamma_2$ imply that 
\[ 
\gamma_1(s, \pi \times \rho,\psi)  =\gamma_2(s, \pi \times \rho,\psi) 
\] 
as desired. This completes the proof of uniqueness and hence the proof of the theorem. 
\vskip 5pt

\vskip 10pt

\section{\bf Adjoint $\gamma$-factors}
Following the above lines, one can also construct a theory of adjoint $\gamma$-factors for $G_2$.  More precisely, for $\pi \in \Irr(G_2)$ and $\chi$ a character of $\GL_1$, we set:
\[  \gamma(s, \pi , {\rm Ad} \times \chi,\psi) := \gamma(s, {\rm Lif}(\pi), \wedge^2 \times \chi, \psi) / \gamma(s, {\rm Lif}(\pi) \times \chi, \psi)\]
where the first term on the RHS refers to the twisted exterior square L-factors of $\GL_7$ and the second term is the twisted standard $\gamma$-factor. 
Using Proposition \ref{P:SC}, one sees that this system of adjoint $\gamma$-factors satisfies the property of multiplicativity and is the right one for unramified representations $\pi$.
It also satisfies and   is  uniquely characterized by the analogous properties in Theorem \ref{T:gamma}. In particular, in the situations of Theorem \ref{T:gamma}(ii),  the global partial adjoint L-function has meromorphic continuation and satisfies the analogous global functional equation with these local $\gamma$-factors. 
\vskip 10pt

Using the fact that the LLC for $\GL_r$ respects twisted exterior square $\gamma$-factors, a result of Henniart \cite{He}, one has: 

\vskip 5pt

\begin{thm}
The twisted adjoint $\gamma$-factor defined above is compatible with the LLC, i.e for $\pi \in \Irr(G_2)$ and $\chi$ a character of $\GL_1$,  one has:
 \[  \gamma(s, \pi, {\rm Ad} \times \chi) = \gamma(s, \wedge^2( {\rm std} \circ \phi_{\pi}) \otimes \chi, \psi) / \gamma(s, ({\rm std} \circ \phi_{\pi} )\otimes \chi, \psi), \]
where $\phi_{\pi}: WD_F \longrightarrow G_2(\C)$ and $\phi_{\tau}: WD_F \longrightarrow \GL_r(\C)$  are the L-parameters of $\pi$ and $\tau$ respectively and ${\rm std}: G_2(\C) \longrightarrow \GL_7(\C)$ is the standard degree 7 irreducible  representation of $G_2(\C)$. 
\end{thm}

\end{document}